\documentclass[11pt]{amsart}

\usepackage{amsmath}
\usepackage{amssymb}
\usepackage[all]{xy}
\usepackage{sty1}

\makeatletter
\@addtoreset{equation}{section}
\makeatother

\def\bloneg{\mathrm{L}^1(G)}

\def\blonek{\mathrm{L}^1(K)}

\def\falg{\mathrm{A}(G)}
\def\falgg{\mathrm{A}(G\cross G)}
\def\fal#1{\mathrm{A}(#1)}
\def\falggd{\mathrm{A}(G\cross G\!:\!\Delta)}
\def\fdelg{\mathrm{A}_{\Del}(G)}
\def\fdelh{\mathrm{A}_{\Del}(H)}
\def\fdel#1{\mathrm{A}_{\Del}(#1)}
\def\fdelgg{\mathrm{A}_{\Del}(G\cross G)}
\def\fdelggs{\mathrm{A}_{\Del}(G\times G)}
\def\fdelgh{\mathrm{A}_{\Del}(G\cross H)}
\def\fdelsqg{\mathrm{A}_{\Del^2}(G)}
\def\fdelng{\mathrm{A}_{\Del^n}(G)}
\def\fdelngg{\mathrm{A}_{\Del^n}(G\cross G)}
\def\fdelbg#1{\mathrm{A}_{\Del^{#1}}(G)}
\def\fcdelg{\mathrm{A}_{\check{\Del}}(G)}
\def\fgamg{\mathrm{A}_{\gam}(G)}
\def\fsalg{\mathrm{B}(G)}
\def\id{\mathrm{id}}
\def\ideal{\mathrm{I}}
\def\ran{\mathrm{ran}}
\def\t{\mathrm{t}}
\def\vng{\mathrm{VN}(G)}

\def\fnorm#1{\left\|#1\right\|_{\mathrm{A}}}
\def\fdelnorm#1{\left\|#1\right\|_{\mathrm{A}_\Del}}
\def\fdelsqnorm#1{\left\|#1\right\|_{\mathrm{A}_{\Del^2}}}
\def\fdelnnorm#1{\left\|#1\right\|_{\mathrm{A}_{\Del^n}}}
\def\fgamnorm#1{\left\|#1\right\|_{\mathrm{A}_{\gam}}}

\begin{document}

\newtheorem{spectrum}{Proposition}[section]
\newtheorem{optensprod2}[spectrum]{Proposition}
\newtheorem{basicprop}[spectrum]{Proposition}
\newtheorem{idealr}[spectrum]{Theorem}
\newtheorem{idealr1}[spectrum]{Corollary}

\newtheorem{tnormcoeff}{Lemma}[section]
\newtheorem{fdelgchar}[tnormcoeff]{Theorem}
\newtheorem{fdelgchar1}[tnormcoeff]{Corollary}
\newtheorem{fdelgchar2}[tnormcoeff]{Corollary}
\newtheorem{optensprod}[tnormcoeff]{Proposition}
\newtheorem{translation}[tnormcoeff]{Proposition}
\newtheorem{fdelngchar}[tnormcoeff]{Theorem}

\newtheorem{comweakamen}{Theorem}[section]
\newtheorem{comweakamen2}[comweakamen]{Corollary}
\newtheorem{restricth}[comweakamen]{Lemma}
\newtheorem{fident}[comweakamen]{Lemma}
\newtheorem{weakamen}[comweakamen]{Theorem}
\newtheorem{weird}[comweakamen]{Lemma}
\newtheorem{opensub}[comweakamen]{Lemma}
\newtheorem{weakamens}[comweakamen]{Theorem}
\newtheorem{weakamens2}[comweakamen]{Corollary}
\newtheorem{opamen}[comweakamen]{Theorem}
\newtheorem{specsynth}[comweakamen]{Theorem}

\newtheorem{convolution}{Theorem}[section]

\title[Convolutions and coset spaces]
{Convolutions on compact groups and
Fourier algebras of coset spaces}

\author{Brian E. Forrest, Ebrahim Samei and Nico Spronk}

\begin{abstract}
In this note we study two related questions.  (1) For a compact group
$G$, what are the ranges of the convolution maps on $\falgg$
given for $u,v\iin\falg$ by $u\cross v\mapsto u\con\check{v}$
($\check{v}(s)=v(s^{-1})$) and $u\cross v\mapsto u\con v$?
(2) For a locally compact group $G$ and a compact subgroup
$K$, what are the amenability properties of the Fourier
algebra of the coset space $\fal{G/K}$?  The algebra $\fal{G/K}$
was defined and studied by the first named author.

In answering the first question, we obtain for compact groups which do not
admit an abelian subgroup of finite index, some new subalgebras of $\falg$.
Using those algebras we can find many instances in
which $\fal{G/K}$ fails the most rudimentary amenability property:
operator weak amenability.  However, using different techniques, we show
that if the connected component of the identity of $G$ is abelian, then
$\fal{G/K}$ always satisfies the stronger property that it is hyper-Tauberian,
which is a concept developed by the second named author.
We also establish a criterion which characterises operator amenability
of $\fal{G/K}$ for a class of groups which includes the maximally almost 
periodic groups.  Underlying our calculations are some refined techniques
for studying spectral synthesis properties of sets for Fourier algebras.  We even
find new sets of synthesis and nonsynthesis for Fourier algebras of some
classes of groups.
\end{abstract}

\maketitle

\footnote{{\it Date}: \today.


Research of the the first named author supported by NSERC Grant 90749-00.
Research of the second named author supported by an NSERC Post Doctoral 
Fellowship.  
Research of the the third named author supported by NSERC Grant 312515-05.}

In a recent article \cite[\S 3.3]{spronk2}, the third named author posed a 
question which, for a compact group $G$, reduces to asking if the maps
from $\falgg$ to $\falg$, given on elementary functions by
$u\cross v\mapsto u\con\check{v}$ and $u\cross v\mapsto u\con v$,
are surjective.  We show in Sections \ref{sec:fdelg} and \ref{sec:fdelcg}
that this is not the case when $G$ does not admit an abelian subgroup
of finite index.  Moreover, the ranges of both maps are quite different:
the first gives us a new algebra $\fdelg$, and the second gives us
an algebra $\fgamg$, which was originally discovered by B.\ E.\ Johnson 
\cite{johnson}.  It is worth noting that Johnson used $\fgamg$ in a very 
clever way to show that for compact groups $G$,
$\falg$ is generally not amenable, in fact not even weakly amenable.
Johnson's results were surprising when his article was published, since 
at the time
the expectation was that amenability of $\falg$ would be equivalent to 
classical amenability for the underlying group $G$. Recent applications of the 
theory of operator spaces to the study of $\falg$ have given us a much 
better understanding of why $\falg$ fails to be amenable as a Banach algebra 
for any group which does not contain abelian subgroups of finite index. 
However, it is one of the main themes of this paper that subalgebras of 
$\falg$ such as $\fgamg$ can shed further light on the amenability problem and
allow us to deduce much more about the relationship between $\falg$ and $G$.

Motivated by Johnson's beautiful theorem that 
the group algebra $\bloneg$ of a locally compact group $G$ is amenable 
as a Banach algebra if and only if $G$ is amenable, together with the 
spectacular failure of the natural analog of this result for $\falg$, 
Z.-J. Ruan demonstrated the tremendous value in recognising
the operator space structure on $\falg$ by introducing 
the concept of operator amenability and then using this to show that
$\falg$ is operator amenable if and only if $G$ is amenable \cite{ruan}.
This was followed by the third named author \cite{spronk}, and, independantly
the second named author \cite{samei},
each establishing that $\falg$ is always operator weakly amenable.
The general question of when $\falg$ is weakly 
amenable as a Banach algebra remains open.   However, using, in part, techniques
developed in the present article, the authors have recently shown $\falg$ fails
to be weakly amenable for any $G$ which contains a connected nonabelian
compact subgroup.  This will appear in another article, soon.

For a compact subgroup $K$ of $G$, the Fourier algebra
of the coset space, $\fal{G/K}$ was described by the first named
author \cite{forrest}.  $\fal{G/K}$ may be simultaneously viewed as an algebra 
of continuous functions on the coset space $G/K$ and as a sublgebra
of $\falg$. The latter view allows us to define a natural operator space 
structure on $\fal{G/K}$.
It was shown in \cite{forrest} that many properties of $\fal{G/K}$
associated with amenability, such as existence of a bounded approximate
identity, and factorisation, are closely linked to such properties of $\falg$.
Thus we are naturally led to consider the amenability properties
of $\fal{G/K}$.  Surprisingly, positive results are rather sparse,
even in the category of operator spaces. 
In Section \ref{sec:weakamen} we establish that whenever
$G$ has a compact connected nonabelian subgroup $K$, then
there exists a compact subgroup $K^*$ of $G\cross G$
such that $\fal{G\cross G/K^*}$ is not operator weakly amenable.
This contrasts sharply the positive result of \cite{spronk,samei}, 
mentioned above.
As a complement, we establish in Section \ref{ssec:hypertaub}
that $\fal{G/K}$ is hyper-Tauberian
whenever the connected component of the identity 
$G_e$ in $G$ is abelian.  
The hyper-Tauberian
property, developed by the second named author \cite{samei2},
implies weak amenability of a commutative semisimple Banach algebra.
Note that this result does not address the operator space
structure of $\falg$ at all.  This generalises a result of 
the first named author and V.\ Runde \cite{forrestr} 
that $\falg$ is weakly amenable when $G_e$ is abelian.
In Section \ref{ssec:opamen} we obtain, for certain groups which
we call [MAP]$_K$-groups,
a characterisation of when $\fal{G/K}$ is operator
amenable.  Using related techniques, we also find new sets of
local synthesis for $\falgg$ when $G$ has an abelian connected component
of the identity, ones which are not known to be
in the closed coset ring of $G\cross G$.

Many of the results of Section \ref{sec:weakamen} rely heavily
on constructions relating to convolutions on compact groups.
In Section \ref{sec:twisted} we develop a general 
framework in which to view the `twisted'
convolution $f\cross g\mapsto f\con \check{g}$ on functions
on compact $G$.  We realise the image of this map
as a special example of functions on a coset space.
We find that our general framework naturally accommodates
an easy, though far reaching generalisation of a result
relating spectral information between various different algebras,
obtained for abelian compact groups by N. Th.\ Varapolous 
\cite{varopoulos}, and generalised to arbitrary compact groups
by L. Turowska and the third named author; see Theorem
\ref{theo:idealr}.  Since it takes us little extra effort,
we prove our spectral results for a class we call [MAP]$_K$-groups, a class which
includes maximally almost periodic groups.
In Section \ref{sec:fdelg} we apply our twisted
convolution framework to $\falgg$.  In doing so, we obtain
not one, but an infinite sequence of new subalgebras 
of $\falg$, when $G$ does not admit an abelian subgroup of finite 
index.  In Section \ref{sec:fdelcg} we provide a framework
for convolutions on compact $G$, and show how it differs
from the twisted convolution when applied to $\falgg$.
In effect, we have an alternate method to obtain the algebra
$\fgamg$ of Johnson.

\subsection{Background and notation}

The {\it Fourier algebra} $\falg$ of a locally compact group $G$
was defined by Eymard \cite{eymard}.  For compact $G$, there is
an alternative description in \cite[Chap.\ 34]{hewittrII}.
That the two descriptions coincide can be seen by comparing
\cite[p.\ 218]{eymard} with \cite[(34.16)]{hewittrII}.  We note
that the Fourier algebra is closed under both 
group actions of left and right translations
$(s,u)\mapsto s\con u,\,s\mult u:G\cross\falg\to\falg$, given by
\[
s\con u(t)=u(s^{-1}t),\qquad s\mult u(t)=u(ts).
\]
Moreover, these actions are continuous in $G$ and isometric on $\falg$.
We note that $\falg$ admits a von Neumann algebra $\vng$ as its dual space.
As such it comes equipped with a natural operator space structure.
See \cite{effrosrB}, for example, for more on this.  We use the
same definitions as \cite{effrosrB} for {\it completely bounded map},
{\it complete isometry} and 
{\it complete quotient map}.  We note that right and left translations
on $\falg$, being the preadjoints of multiplications by unitaries
on $\vng$, are complete isometries.

Our main references for amenability are \cite{johnsonm} and \cite{rundeb}.
A Banach algebra $\fA$ is said to be {\it amenable} if for any Banach 
$\fA$-bimodule $\fX$, and any bounded derivation $D:\fA\to\fX^*$,
where $\fX^*$ is the dual space with adjoint module actions, $D$
is inner.  For a commutative Banach algebra $\fA$, we say
$\fA$ is {\it weakly amenable} if for any symmetric $\fA$-bimodule
$\fX$, the only bounded derivation $D:\fA\to\fX$ is $0$; this
is equivalent to having the same happen for $\fX=\fA^*$.  Weak
amenability for commutative Banach algebras was introduced
in \cite{badecd}.  For both amenability and weak amenability
there are some homological characterisations; see 
\cite{helemski,curtisl,groenbaek,rundeb}.

Operator space notions of amenability and weak amenability
were introduced in \cite{ruan} and in \cite{forrestw} respectively, 
specifically for use with $\falg$.
If $\fA$ is a commutative Banach algebra which is also an operator space
we say $\fA$ is a {\it completely contractive Banach algebra}
if the multiplication map $\fA\cross\fA\to\fA$ is completely
contractive in the sense of \cite[Chap.\ 7]{effrosrB}.
An operator space $\fV$ is a {\it completely contractive $\fA$-module}
if it is an $\fA$ module for which
the module maps $\fA\cross\fV,\fV\cross\fA\to\fV$ are completely
contractive.  The class of completely contractive $\fA$-modules
is closed under taking dual spaces with adjoint actions.
We say $\fA$ is {\it operator (weakly) amenable} if every
completely bounded derivation $D:\fA\to\fV^*$ (with $\fV=\fA$)
is inner (zero).  Many of the homological characterisations
alluded to above, carry over to this setting, though with Banach space
projective tensor products $\otimes^\gam$ replaced by
operator space projective tensor products $\what{\otimes}$.

Let $\fA$ be a commutative semisimple (completely contractive) Banach algebra.
Suppose $\fA$ is regular on its spectrum $X=\Sig_\fA$;
we regard $\fA$ as an algebra of functions on $X$.
If $\vphi\in\fA^*$ we define
\[
\supp{\vphi}=\left\{x\in X:
\begin{matrix} \text{for every neighbourhood }U\oof x 
\text{ there is }f \\
\iin\fA\text{ such that }\supp{f}\subset U\aand \vphi(f)\not=0
\end{matrix}\right\}.
\]
Here $\supp{f}=\wbar{\{x\in X:f(x)\not=0\}}$.  An operator
$T:\fA\to\fA^*$ is called a {\it local map} if 
\[
\supp{Tf}\subset\supp{f}
\]
for every $f\iin\fA$.  We say $\fA$ is {\it (operator) hyper-Tauberian} if 
every (completely bounded) bounded local map $T:\fA\to\fA^*$ is an 
$\fA$-module map.  This concept
was developed in \cite{samei2} to study
the reflexivity of the (completely bounded) derivation space of 
$\fA$, and it generalises (operator) weak amenability.

For a commutative semisimple Banach algebra which is regular on its spectrum
we have the following implications.
\[
\xymatrix{
\text{amenable } \ar@{=>}[r]
& \text{ weakly amenable } 
& \text{ hyper-Tauberian} \ar@{=>}[l] 
}
\]
Moreover, if $\fA$ is a completely contractive Banach algebra, each
property implies its operator analogue, and the operator analogues
satisfy the same implications.

\section{Algebras on coset spaces and a twisted
convolution on compact groups}\label{sec:twisted}

\subsection{The basic construction}\label{ssec:basic}
Let $G$ be a locally compact group.
Let $\fA(G)$ be a Banach
algebra of continuous functions on $G$ which
is closed under right translations and such that
for any $f\iin\fA(G)$ we have
\begin{align*}
&\bullet\;\norm{s\mult f}=\norm{f}\text{ for any }s\iin 
G,\aand \\
&\bullet\; s\mapsto s\cdot f\text{ is continuous}.
\end{align*}
If, moreover, $\fA(G)$ is an operator space, we
want that for each $t\iin G$ that $f\mapsto t\mult f$
is a complete isometry.

If $K$ is a compact subgroup of $G$ we let
\[
\fA(G\!:\!K)=\{f\in\fA(G):k\mult f=f\text{ for each }k\iin K\}
\]
which is a closed subalgebra of $\fA(G)$ whose elements
are constant on left cosets of $K$.  We let $G/K$
denote the space of left cosets with the quotient topology.
We define two maps
\begin{align*}
P:\fA(G)\to\fA(G),\qquad & Pf=\int_K k\mult f dk \aand \\
M:\fA(G\!:\! K)\to\fC_b(G/K),\qquad & Mf(sK)=f(s).
\end{align*}
The map $P$ is to be regarded as a Bochner integral
over the normalised Haar measure on $K$; its
 range is $\fA(G\!:\!K)$ and $P$ is a (completely) contractive
 projection.
The map $M$ is well-defined by comments above,
and its range consists of continuous functions since
$\fA(G\!:\!K)\subset\fC_b(G\!:\!K)$.
We note that $M$ is an injective homomorphism
and denote its range by
$\fA(G/K)$.  We assign a norm (operator space structure)
to $\fA(G/K)$ in such a way that $M$ is a (complete)
isometry.  We finally define two maps
\[
N=M^{-1}:\fA(G/K)\to\fA(G)\aand
\Gam=M\comp P:\fA(G)\to\fA(G/K)
\]
so $N$ is a (completely) isometric homomorphism and
$\Gam$ is a (complete) quotient map.

Let us record some basic properties of $\fA(G/K)$.
For a commutative Banach algebra $\fA$, we let $\Sigma_\fA$
denote its Gelfand spectrum.

\begin{spectrum}\label{prop:spectrum}
{\bf (i)} Suppose $\fA(G)$ 
is regular on $G$ and $G$ separates
the points of $\fA(G)$.  Then $\fA(G/K)$ is regular
on $G/K$ and $G/K$ separates the points of $\fA(G/K)$.  

Moreover, if $\Sig_{\fA(G)}\cong G$ via evaluation maps, and
$K$ is a set of spectral synthesis for
$\fA(G)$, then $\Sigma_{\fA(G/K)}\cong G/K$ via evaluation maps.

{\bf (ii)} If the subalgebra $\fA_c(G)$ of compactly supported
elements of $\fA(G)$ is dense in $\fA(G)$,
then the algebra $\fA_c(G/K)$ of
compactly supported elements in $\fA(G/K)$ is dense in $\fA(G/K)$.
\end{spectrum}

We remark that (i) applies to the Fourier algebra $\falg$
for any locally compact group $G$ and compact subgroup $K$
by \cite{herz} or \cite{takesakit}.

\medskip
\proof {\bf (i)} Same as \cite[Thm.\ 4.1]{forrest}.  

{\bf (ii)} It is obvious that $\Gam\fA_c(G)=\fA_c(G/K)$.
If $w\iin\fA(G)$ is the limit
of a sequence $(w_n)\subset\fA_c(G)$, then 
$\norm{\Gam w-\Gam w_n}\leq \norm{w-w_n}\longrightarrow 0$
as $n\to\infty$.
\endpf

For Fourier algebras, we have the important identification
\[
\falg\what{\otimes}\fal{H}\cong\fal{G\cross H}
\]
via $f\otimes g\mapsto f\cross g$, where $\what{\otimes}$
denotes the operator projective tensor product.  See
\cite{effrosr2}.  This is known to fail when the usual
projective tensor product $\otimes^\gam$ is used \cite{losert}.

\begin{optensprod2}\label{prop:optensprod2}
If $G$ and $H$ are locally compact groups with respective compact 
subgroups $K$ and $L$, then
there is a complete isometry identifying
\[
\fal{G/K}\what{\otimes}\fal{H/L}\cong\fal{G\cross H/K\cross L}
\]
given on elementary tensors by $f\otimes g\mapsto f\cross g$.
\end{optensprod2}

Note that there is a natural homeomorphism $G\cross H/K\cross L
\cong G/K\cross H/L$, which is in fact a $G\cross H$-space morphism.

\medskip
\proof We identify $\fal{G/K}\cong\fal{G\!:\! K}$, etc.
We have the following commuting digram
\[
\xymatrix{
\falg\what{\otimes}\fal{H} \ar[rr]^{f\otimes g\mapsto f\times g}
\ar@<1ex>[d]^{P_K\otimes P_L}
& & \fal{G\cross H} \ar[d]^{P_{K\times L}} \\
\fal{G\!:\!K}\what{\otimes}\fal{G\!:\!K} \ar@<1ex>@{^{(}->}[u]
\ar@{-->}[rr]^{f\otimes g\mapsto f\times g}
& & \fal{G\cross H\!:\! K\cross L}}
\]
where the inclusion map $\fal{G\!:\!K}\what{\otimes}\fal{G\!:\!K}
\hookrightarrow\falg\what{\otimes}\fal{H}$ is a complete isometry,
since each of $P_K$ and $P_L$ are complete quotient projections.
Notice that this inclusion map is a right inverse to $P_K\otimes P_L$.  
The desired map 
from $\fal{G\!:\!K}\what{\otimes}\fal{G\!:\!K}$ to
$\fal{G\cross H\!:\! K\cross L}$ completes this digram, and may be
realised as the composition of the inclusion map (up arrow), 
with the identification map (top arrow), and then $P_{K\times L}$.
It is clear that $P_{K\times L}|_{\wbar{\spn}\fal{G:K}\times\fal{H:L}}=\id$.
Hence the desired map is injective. 
Since the map identifying $\falg\what{\otimes}\fal{H}\cong\fal{G\cross H}$
is a complete isometry, and $P_K\otimes P_L$ and $P_{K\times L}$
are complete quotient maps, our desired map is a complete quotient map.
Hence we have an injective complete quotient map, i.e.\ a complete isometry.
\endpf


\subsection{A twisted convolution}\label{ssec:twisted}
Let $G$ be a {\it compact} group for the remainder of this section. 
We use $G\cross G$ in place of $G$ above, and
$K=\Del=\{(s,s):s\in G\}$.  Since the map
\begin{equation}\label{eq:ggmoddel}
(s,e)\Del\mapsto s:(G\cross G)/\Del\to G 
\end{equation}
is a homeomorphism,
we identify the coset space with $G$.  
We observe that in this case that the map $P:\fA(G\cross G)
\to\fA(G\cross G)$ satisfies
\[
Pw(s,t)=\int_G w(sr,tr)dr=\int_G w(st^{-1}r,r)dr
\]
and the maps
$M:\fA(G\cross G\!:\!\Del)\to\fC(G)$ and $N:\fC(G)\to\fC(G\cross G\!:\!\Del)$
satisfy
\begin{equation}\label{eq:nandm}
Mw(s)=w(s,e)\quad\aand \quad Nf(s,t)=f(st^{-1}).
\end{equation}
We denote the range of $M$ by $\fA_\Del(G)$, and then
$N|_{\fdelg}=M^{-1}$.
The map $\Gam=M\comp P$, above, can
be regarded as a `twisted' convolution, for if $\fA(G\cross G)$
contains an elementary function $f\cross g$, then for $s\iin G$
\[
\Gam(f\cross g)(s)=\int_G f\cross g(st,t)dt
=\int_Gf(st)g(t)dt=f\con\check{g}(s).
\]

We list some examples of $\fA(G\cross G)$ and $\fA_\Del(G)$.

\medskip
{\bf (i)}  If $\fA(G\cross G)=\fC(G\cross G)$, then
$\fC_\Del(G)=\fC(G)$, by easy computation.

{\bf (ii)} Let $\mathrm{V}(G\cross G)=\fC(G)\otimes^h\fC(G)$ (Haagerup
tensor product). 
Then $\mathrm{V}_\Del(G)=\falg$, completely 
isometrically.  See \cite{spronkt}.
Note that by Grothedieck's inequality
$\mathrm{V}(G\cross G)=\fC(G)\otimes^\gam\fC(G)$ (projective tensor
product) isomorphically, though not isometrically.  

{\bf (iii)} Let $\fA(G\cross G)=\falg\otimes^\gam\falg$.
Then $\fA_\Del(G)=\mathrm{A}_\gam(G)$ is a subalgebra of $\falg$
considered by B. Johnson \cite{johnson}.  He used it to study the
amenability of $\falg$. 

{\bf (iv)} Consider $\falgg\cong\falg\what{\otimes}\falg$ (operator
projective tensor product).  
The algebra $\fdelg$, thus formed,
will be an essential object of our study. 


\medskip
We summarise some basic properties of the algebras $\fA_\Del(G)$
which clearly apply to the examples above.  We say a norm (operator space
structure) $\alp$, on $\fX\otimes\fY$ is (operator)  homogeneous
if for every pair of (completely) contractive linear maps $S:\fX\to\fX'$,
$T:\fY\to\fY'$ we have that $S\otimes T$ extends to a (complete)
contraction from $\fX\otimes^\alp\fY$ to $\fX'\otimes^\alp\fY'$.

\begin{basicprop}\label{prop:basicprop}
{\bf (i)} If there is a (completely contractive)
Banach algebra $\fB(G)$ of continuous functions on $G$ 
and a homogeneous (operator) norm (operator space structure)
$\alp$ on $\fB(G)\otimes\fB(G)$ so 
$\fA(G\cross G)=\fB(G)\otimes^\alp\fB(G)$, then
$\fA_\Del(G)$ is a (completely contractive) subalgebra of $\fB(G)$.

{\bf (ii)} If $\fA(G\cross G)$ is closed under left translations
and the translation maps are
continuous in $G\cross G$ and (completely) isometric on $\fA(G\cross G)$,
then $\fA_\Del(G)$ is closed under both left and right translations,
and the translations are continuous on $G$, and (completely) isometric
on $\fA_\Del(G)$.
\end{basicprop}

\proof {\bf (i)} The map $M:\fA(G\cross G\!:\!\Del)\to\fA_\Del(G)$ is the 
restriction of the
slice map $\id\otimes\del_e:\fA(G)\otimes^\alp\fA(G)\to\fA(G)$
where $\del_e$ is the evaluation functional at $e$.  

{\bf (ii)} For $f\in\fA_\Del(G)$,
$r\iin G$ and any $(s,t)\in G\cross G$ we have
\[
N(r\mult f)(s,t)=f(st^{-1}r)=f(s(r^{-1}t)^{-1})=(e,r)\con Nf.
\]
Our assumptions assure that
the space $\fA(G\cross G\!:\!\Del)$ is closed under left translations
and each translation map is a (complete) isometry.  
For left translations  we note that for $r\iin G$ and
$f\in\fdelg$, we have
\[
N(r\con f)=(r,e)\con Nf
\]
and we argue as above.  \endpf

\subsection{Relationships between ideals}\label{ssec:idealr}
Let $\fA$ be a Banach algebra contained in $\fC_0(X)$
for some locally compact Hausdorff space $X$.  We define
for any closed subset $E$ of $X$
\begin{align*}
\ideal_\fA(E)&=\{f\in\fA:f(x)=0\ffor x\in E\},  \\
\ideal_\fA^0(E)&=\{f\in\fA:\supp{f}\cap E=\varnothing\
\aand\supp{f}\text{ is compact}\}, \aand \\
\mathrm{J}^0_\fA(E)&=\{f\in\ideal_\fA(E):\supp{f}\text{ is compact}\}
\end{align*}
where $\supp{f}=\wbar{\{x\in X:f(x)\not=0\}}$.
If $X=\Sig_\fA$, and $\fA$ is regular on $X$
we say $E$ is a set of {\it spectral synthesis}
for $\fA$ if $\wbar{\ideal_\fA^0(E)}=\ideal_\fA(E)$, and of
{\it local synthesis} if $\wbar{\ideal_\fA^0(E)}=\wbar{\mathrm{J}^0_\fA(E)}$.

If $G$ is a compact group
we let $\what{G}$ denote the dual object of $G$, 
a set of representatives, one from each unitary equivalence class, 
of irreducible continuous
representation of $G$. If $\pi\in\what{G}$, we fix an orthonormal basis
$\{\xi_1^\pi,\dots,\xi_{d_\pi}^\pi\}$ for $\fH_\pi$ and define
\begin{equation}\label{eq:piij}
\pi_{ij}:G\to\Cee,\quad \pi_{ij}(s)=\inprod{\pi(s)\xi^\pi_j}{\xi^\pi_i}
\end{equation}
for $i,j=1\dots d_\pi$.  We recall the well-known fact that
\begin{equation}\label{eq:trigfunc}
\fT(G)=\spn\{\pi_{ij}:\pi\in\what{G},i,j=1,\dots,d_\pi\}
\end{equation}
is uniformly dense in $\fC(G)$.

If $G$ is not compact, we let $\what{G}_f$ denote the finite dimensional
part of the dual object.  We let $\fT_f(G)$ be defined as the span
of matrix coefficients of $\what{G}_f$, analogously as above.
The almost periodic compactification is the compact group
\[
G^{ap}=\wbar{\left\{(\pi(s))_{\pi\in\what{G}_f}:s\in G\right\}}
\subset\prod_{\pi\in\what{G}_f}\wbar{\pi(G)}.
\]
There is a canonical identification $\fT(G^{ap})\cong\fT_f(G)$.
We say $G$ is {\it maximally almost periodic} [MAP] if $\what{G}_f$,
or equivalently $\fT_f(G)$, separates the points of $G$.
If $K$ is a compact subgroup of $G$, we say $G$ is [MAP]$_K$
if the map $k\mapsto(\pi(k))_{\pi\in\what{G}_f}:K\to G^{ap}$
is injective.  Clearly, a [MAP] group is [MAP]$_K$ for any
compact subgroup $K$.

The following is an adaptation of \cite[Thm.\ 3.1]{spronkt}.
A change in perspective allows us to gain not only
more general, but finer results than in \cite{spronkt}.

\begin{idealr}\label{theo:idealr}
Let $G$ be a locally compact group, $K$ a compact subgroup 
so $G$ is [MAP]$_K$, and
$\fA(G)$ be as in Section \ref{ssec:basic} and additionally
satisfy that $\fT_f(G)\fA(G)\subset\fA(G)$.  
If $E$ is a closed subset
of $G/K$ let $E^*=\{s\in G:sK\in E\}$.  Then

{\bf (i)} $\Gam\ideal_{\fA(G)}(E^*)=\ideal_{\fA(G/K)}(E)$, and

{\bf (ii)} $\ideal_{\fA(G)}(E^*)$ is the closed ideal generated
by $N\ideal_{\fA(G/K)}(E)$.
\end{idealr}

Note that if $G$ is compact then
in the case of Section \ref{ssec:twisted}
we have $E^*=\{(s,t)\in G\cross G:st^{-1}\in E\}$
via the identification (\ref{eq:ggmoddel}).

\medskip
\proof  We will let $\ideal(E^*)=\ideal_{\fA(G)}(E^*)$
and $\ideal(E)=\ideal_{\fA(G/K)}(E)$, below.

{\bf (i)}  It is clear that
\[
\Gam\ideal(E^*)\subset\ideal(E)\aand
N\ideal(E)\subset \ideal(E^*).
\]
Thus
\[
\ideal(E)=\Gam\comp N \ideal(E)
\subset \Gam\ideal(E^*)\subset \ideal(E).
\]

{\bf (ii)}  Let $w\in\ideal(E^*)$.  For each 
$\pi\iin \what{G}$
we define `matrix-valued' functions
$w^\pi,\til{w}^\pi:G\to\fB(\fH_\pi)$ by
\[
w^\pi(s)=\int_K w(sk)\pi(k)dk\quad\aand\quad
\til{w}^\pi(s)=\pi(s)w^\pi(s).
\]
Then for any $i,j=1,\dots,d_\pi$ we let $w_{ij}^\pi=
\inprod{w^\pi(\cdot)\xi^\pi_j}{\xi^\pi_i}$ and we have
\[
w^\pi_{ij}=\pi_{ij}|_K\mult w
\]
where $f\mult w=\int_G f(k)k\mult wdk$ for any $f\iin\blonek$.
We note that since $w\in\ideal(E^*)$, $k\mult w\in
\ideal(E^*)$ for any $k\iin K$ and hence
$f\mult w\in\ideal(E^*)$ for any $f\iin\blonek$.
Thus each $w^\pi_{ij}\in\ideal(E^*)$.
Now for any $s\iin G$, $i,j=1,\dots,d_\pi$ we have
\[
\til{w}^\pi_{ij}(s)=\sum_{l=1}^{d_\pi}\pi_{il}(s)w^\pi_{lj}(s),\quad
\text{i.e. }\til{w}^\pi_{ij}=\sum_{k=1}^{d_\pi}\pi_{ik}w^\pi_{lj}.
\]
Hence since $\fT_f(G)\fA(G)\subset\fA(G)$ we have that
each $\til{w}^\pi_{ij}\in\ideal(E^*)$.  However,
it is easily seen that $\til{w}^\pi(sk)=\til{w}^\pi(s)$ for any
$s\iin G$ and $k\iin K$, so each $\til{w}^\pi_{ij}\in\fA(G\!:\!K)$.
Thus each
\[
\til{w}^\pi_{ij}\in\ideal(E^*)\cap\fA(G\!:\!K)
=N\ideal(E).
\]
Now we use the relation $w^\pi(s)=\pi(s^{-1})\til{w}^\pi(s)$,
for $s\iin G$, to see that for each $i,j$ we have
\[
w^\pi_{ij}=\sum_{l=1}^{d_\pi}\check{\pi}_{il}\til{w}^\pi_{lj}
=\sum_{l=1}^{d_\pi}\bar{\pi}_{li}\til{w}^\pi_{lj}
\]
so $w^\pi_{ij}=\pi_{ij}|_K\mult w$ lies in the ideal generated
by $N\ideal(E)$.

Since $G$ is [MAP]$_K$, we may regard $K$ as a closed subgroup of $G^{ap}$.
We have that $\fT_f(G)|_K=\fT(G^{ap})|_K$ is a conjugation-closed
point-separating subalgebra of $\fC(K)$, thus
is uniformly dense in $\fC(K)$, and hence 
norm dense in $\blonek$.   
Thus there is a bounded approximate identity $(f_\beta)$
for $\blonek$ for which each $f_\beta\in\fT_f(G)|_K$.
Then for each $\beta $ we have
\[
f_\beta\mult w\in\spn\{\pi_{ij}|_K\mult w:\pi\in\what{G}
\aand i,j=1,\dots,d_\pi\}
\]
and is thus in the ideal generated by $N\ideal(E)$.
Hence, since $\fA(G)$ is an essential $\blonek$-module, we find that
$w=\lim_\beta f_\beta \mult w$ and thus is in the closed ideal generated by 
$N\ideal(E)$.  \endpf

For any subalgebra $\fB$ of a commutative normed algebra $\fA$, we let
\[
\langle \fB\rangle=\wbar{\spn}\{ab_1+b_2:a\in\fA, b_1,b_2\in\fB\}
\aand 
\fB^2=\spn\{ab:a,b\in\fB\}.
\]

\begin{idealr1}\label{cor:idealr1}
With $G$, $K$, $\fA(G)$, $E\subset G/K$ and $E^*\subset G$
as in Theorem \ref{theo:idealr}, we have that

{\bf (i)} $\ideal_{\fA(G)}(E^*)$ has a bounded approximate
identity (b.a.i.) if and only if $\ideal_{\fA(G/K)}(E)$ has
a b.a.i.;

{\bf (ii)} $\wbar{\ideal^0_{\fA(G)}(E^*)}=\ideal_{\fA(G)}(E^*)$
if and only if $\wbar{\ideal_{\fA(G/K)}^0(E)}=\ideal_{\fA(G/K)}(E)$; and

{\bf (ii')} $\wbar{\ideal^0_{\fA(G)}(E^*)}=\wbar{\mathrm{J}^0_{\fA(G)}(E^*)}$
if and only if 
$\wbar{\ideal_{\fA(G/K)}^0(E)}=\wbar{\mathrm{J}^0_{\fA(G/K)}(E)}$.

{\bf (iii)} 
$\wbar{\ideal_{\fA(G)}(E^*)^2}=\ideal_{\fA(G)}(E^*)$
if and only if $\wbar{\ideal_{\fA(G/K)}(E)^2}=\ideal_{\fA(G/K)}(E)$.
\end{idealr1}

\proof   We let $\ideal(E^*)=\ideal_{\fA(G)}(E^*)$,
$\ideal(E)=\ideal_{\fA(G/K)}(E)$, etc.

{\bf (i)} If $(f_\alp)$ is a b.a.i.\ for $\ideal(E)$,
then $(Nf_\alp)$ is a b.a.i.\ for the subalgebra $N\ideal(E)$.
It is readily checked that $(Nf_\alp)$ is a b.a.i.\ for 
$\langle N\ideal(E)\rangle$.

If $(w_\alp)$ is a b.a.i.\ for $\ideal(E^*)$, then $(\Gam w_\alp)$
is a b.a.i.\ for $\ideal(E)$.  Indeed we have that
$\Gam\comp N=\id$ and, since $P$ is an idempotent,
we have the following `expectation property':
\begin{equation}\label{eq:expectation}
\Gam(w\,Nf)=\Gam(w)f
\end{equation}
for $w\in\fA(G\cross G)$ and $f\in\fA(G/K)$.  Thus if
$f\iin\ideal(E)$ then
\[
\Gam(w_\alp)f-f=\Gam(w_\alp\, Nf-Nf)\longrightarrow 0.
\]

{\bf (ii) \& (ii')}  It is clear that 
\[
f\in\ideal^0(E)\iff Nf\in\ideal^0(E^*)
\aand f\in\mathrm{J}^0(E)\iff Nf\in\mathrm{J}^0(E^*).
\]
It can then be proved exactly as in Theorem \ref{theo:idealr} that
\[
\left\langle N\ideal^0(E)\right\rangle
=\wbar{\ideal^0(E^*)}\aand
\left\langle N\mathrm{J}^0(E)\right\rangle
=\wbar{\mathrm{J}^0(E^*)}.
\]
Indeed, it is sufficient to note that if $\supp{f}$ is compact
then so too is $\supp{(\pi_{ij}|_K\mult f)}$, for each $\pi\in\what{G}_f$
and $i,j=1,\dots,d_\pi$.

{\bf (iii)}  If $\wbar{\ideal(E)^2}=\ideal(E)$,
then $N\bigl(\ideal(E)^2\bigr)
=\bigl(N\ideal(E)\bigr)^2$ is dense in $N\ideal(E)$.
Hence we have
\[
\ideal(E^*)=\langle N\ideal(E)\rangle
=\left\langle (N\ideal(E)\bigr)^2\right\rangle 
=\wbar{\langle N\ideal(E)\rangle^2}
=\wbar{\ideal(E^*)^2}.
\]
Conversely, if $\wbar{\ideal(E^*)^2}=\ideal(E^*)$
we have from the theorem above that
\[
\ideal(E^*)=\wbar{\ideal(E^*)^2}
=\wbar{\ideal(E^*)N\ideal(E)}.
\]
But it follows from (\ref{eq:expectation}) that
\[
\ideal(E)=\Gam\ideal(E^*)
=\wbar{\Gam\bigl(\ideal(E^*)\,N\ideal(E)\bigr)} 
=\wbar{\Gam\ideal(E^*)\,\,\ideal(E)}
=\wbar{\ideal(E)^2}\subset\ideal(E)
\]
whence $\ideal(E)=\wbar{\ideal(E)^2}$. \endpf

\section{Some subalgebras of Fourier algebras}\label{sec:fdelg}

\subsection{The algebra $\fdelg$}
In this section, we will always let $G$ denote a compact group.  
We have
the following characterisation of the Fourier algebra in
\cite[(34.4)]{hewittrII}: for $f\in\fC(G)$
\begin{equation}\label{eq:falg}
f\in\falg\quad\iff\quad
\fnorm{f}=\sum_{\pi\in\what{G}}d_\pi\norm{\hat{f}(\pi)}_1<\infty
\end{equation}
where $\hat{f}(\pi)=\int_Gf(s)\bar{\pi}(s)ds$ and $\norm{\cdot}_1$ denotes
the trace class norm.  We also recall the following orthogonality
relations \cite[(27.19)]{hewittrII}:  
if $\pi,\sig\in\what{G}$ then in the notation of (\ref{eq:piij})
we have
\begin{equation}\label{eq:orthrel}
\int_G\pi_{ij}(s)\wbar{\sig_{kl}(s)}ds
=\frac{1}{d_\pi}\del_{\pi\sig}\del_{jl}\del_{ik}
\end{equation}
where $i,j=1,\dots,d_\pi$, $k,l=1,\dots,d_\sig$ and each $\del_{\alp\beta}$
is the Kronecker $\del$-symbol.

We first wish to characterise $\fdelg$, as defined in the 
previous section.  We will make use of the following lemma.
We let for $\pi,\sig\in\what{G}$, $\pi\times\sig\in\what{G\times G}$
be the Kronecker product representation.  Also,
we let $N:\fC(G)\to\fC(G\cross G)$ be given by $Nf(s,t)=f(st^{-1})$
for $(s,t)\iin G\cross G$, as suggested by (\ref{eq:nandm}).

\begin{tnormcoeff}\label{lem:tnormcoeff}
For any $f\in\fC(G)$ and $\pi\in\what{G}$, we have 
\[
\norm{\what{Nf}(\bar{\pi}\cross\pi)}_p
=\frac{1}{\sqrt{d_\pi}}\norm{\hat{f}(\pi)}_2
\]
where $\norm{\cdot}_p$ is the Schatten $p$-norm for $1\leq p\leq\infty$.
\end{tnormcoeff}

\proof  We first note that
\begin{align}\label{eq:tilmf}
\what{Nf}(\bar{\pi}\cross\pi)
&=\int_G\int_G f(st^{-1})\pi(s)\otimes\bar{\pi}(t)dsdt \notag \\
&=\int_G\int_G f(s)\pi(st)\otimes\bar{\pi}(t)dsdt \notag \\
&=\left[\int_G f(s)\pi(s)\otimes\bar{\pi}(e)ds\right]\comp
\left[\int_G\pi(t)\otimes\bar{\pi}(t)dt\right] \\
&=[\hat{f}(\bar{\pi})\otimes I_{\fH_{\bar{\pi}}}]\comp P_1 \notag
\end{align}
where $P_1$ is a rank $1$ projection, as we shall see below.
Indeed, the Schur orthogonality relations \cite[(27.30)]{hewittrII}
tell us that $\pi\otimes\bar{\pi}$
contains the trivial representation $1$ with multiplicity $1$.  Thus when
we decompose into irreducibles $\pi\otimes\bar{\pi}
\cong\bigoplus_{\sig\in \what{G}}m_\sig\mult\sig$, we obtain
\[
\int_G\pi(t)\otimes\bar{\pi}(t)dt\cong\bigoplus_{\sig\in \what{G}}m_\sig\mult
\int_G\sig(t)dt
\]
where $\int_G\sig(t)dt=0$ unless $\sig=1$, and hence the formula above
reduces to a rank $1$ projection, $P_1$.  Let $\{\xi_1,\dots,\xi_{d_\pi}\}$
be an orthonormal basis for $\fH_{\pi}$; we claim that
\[
P_1=\inprod{\cdot}{\eta}\eta\quad\wwhere\quad\eta=
\frac{1}{\sqrt{d_\pi}}\sum_{k=1}^{d_\pi}\xi_k\otimes\bar{\xi}_k.
\]
Indeed, 
for any $d_\pi\times d_\pi$ unitary matrix $U$ we have
that $\sum_{k=1}^{d_\pi}U\xi_k\otimes\wbar{U\xi}_k
=\sum_{k=1}^{d_\pi}\xi_k\otimes\bar{\xi}_k$.  Hence $\eta$ is a norm $1$
vector, invariant for $\pi\otimes\bar{\pi}$, and the formula for $P_1$ follows.

Thus
\[
\what{Nf}(\bar{\pi}\otimes\pi)=[\hat{f}(\bar{\pi})
\otimes I_{\fH_{\bar{\pi}}}]\comp
P_1=\inprod{\cdot}{\eta}[\hat{f}(\bar{\pi})\otimes I_{\fH_{\bar{\pi}}}]\eta
\]
and, using the standard formula for rank $1$ operators, 
$\norm{\inprod{\cdot}{\zeta}\vartheta}_p=\norm{\zeta}\norm{\vartheta}$,
we obtain
\begin{equation}\label{eq:crossnorm}
\norm{\what{Nf}(\bar{\pi}\otimes\pi)}_p
=\norm{[\hat{f}(\bar{\pi})\otimes I_{\fH_{\bar{\pi}}}]\eta}
_{\fH_{\pi}\otimes^2\fH_{\bar{\pi}}}.
\end{equation}
Letting $\hat{f}(\bar{\pi})_{kl}=\inprod{\hat{f}(\bar{\pi})\xi_l}{\xi_k}$
we have
\[
[\hat{f}(\bar{\pi})\otimes I_{\fH_{\bar{\pi}}}]\eta
=\frac{1}{\sqrt{d_\pi}}\sum_{k=1}^{d_\pi}
\left[\hat{f}(\bar{\pi})\xi_k\right]\otimes\bar{\xi}_k
=\frac{1}{\sqrt{d_\pi}}\sum_{k=1}^{d_\pi}\sum_{l=1}^{d_\pi}
\hat{f}(\bar{\pi})_{kl}\xi_k\otimes\bar{\xi}_l.
\]
Since $\{\xi_l\otimes\bar{\xi}_k:k,l=1\dots,d_\pi\}$ is an orthonormal basis
for $\fH_{\pi}\otimes^2\fH_{\bar{\pi}}$ we obtain that
\begin{equation}\label{eq:finalform}
\norm{[\hat{f}(\bar{\pi})\otimes I_{\fH_{\bar{\pi}}}]\eta}^2
_{\fH_{\pi}\otimes^2\fH_{\bar{\pi}}}
=\frac{1}{d_\pi}\sum_{k=1}^{d_\pi}\sum_{l=1}^{d_\pi}
|\hat{f}(\bar{\pi})_{kl}|^2=\frac{1}{d_\pi}\norm{\hat{f}(\bar{\pi})}_2^2.
\end{equation}
The result is obtained by combining (\ref{eq:crossnorm}) with
(\ref{eq:finalform}) and the fact that $\norm{\hat{f}(\bar{\pi})}_2
=\norm{\hat{f}(\pi)}_2$.  \endpf

We now obtain a characterisation of $\fdelg$ in the spirit of
(\ref{eq:falg}).  

\begin{fdelgchar}\label{theo:fdelgchar}
If $f\in\fC(G)$, then
\[
f\in\fdelg\quad\iff\quad \sum_{\pi\in\what{G}}d_\pi^{3/2}\norm{\hat{f}(\pi)}_2
<\infty.
\]
Moreover, the latter quantity is 
\[
\fdelnorm{f}=
\inf\{\fnorm{w}:w\in\falgg\wwith\Gam w=f\}.
\]
\end{fdelgchar}

\proof Since $N:\fdelg\to\falggd\subset\falgg$ is an isometry we have
\[
f\in\fdelg\quad\iff\quad Nf\in\falgg
\]
in which case $\fdelnorm{f}=\fnorm{Nf}$. 
Recall that $\what{G\cross G}=\{\pi\cross\sig:\pi,\sig\in\what{G}\}$;
see \cite[(27.43)]{hewittrII}, for example.
Note that analogous computations to (\ref{eq:tilmf}),
combined with (\ref{eq:orthrel}), show that
\begin{equation}\label{eq:tilmf1}
\what{Nf}(\pi\cross\sig)=0\iif\sig\not=\bar{\pi}.
\end{equation}
We have that $Nf\in\falgg$ exactly when
\begin{align*}
\fnorm{Nf}&=\sum_{(\sig,\pi)\in\what{G}\times\what{G}}d_\sig d_\pi
\norm{\what{Nf}(\sig\cross\pi)}_1&\text{by (\ref{eq:falg})} \\
&=\sum_{\pi\in\what{G}}d_\pi^2\norm{\what{Nf}(\bar{\pi}\cross\pi)}_1
&\text{by (\ref{eq:tilmf1})} \\
&=\sum_{\pi\in\what{G}}d_\pi^{3/2}\norm{\hat{f}(\pi)}_2&\text{by Lemma
\ref{lem:tnormcoeff}}
\end{align*}
and the latter quantity is finite.  \endpf

We prove some consequences of the result above.  Let us note for
any $d\cross d$ matrix the well-known inequalities
\begin{equation}\label{eq:hstrace}
\frac{1}{\sqrt{d}}\norm{A}_1\leq\norm{A}_2\leq\norm{A}_1.
\end{equation}
These inequalities are sharp with scalar matrices $\alp I$
serving as the extreme case for the left inequality, and rank $1$
matrices serving for the extreme case on the right.

We let $\mathrm{Z}\bloneg$ be the centre of the convolution algebra $\bloneg$.

\begin{fdelgchar1}
We have that $\falg\cap\mathrm{Z}\bloneg=\fdelg\cap\mathrm{Z}\bloneg$ 
with $\fnorm{f}=\fdelnorm{f}$ for $f$ in this space.
\end{fdelgchar1}

\proof  It is well-known that $f\in\mathrm{Z}\bloneg$ only if
$\hat{f}(\pi)=\alp_{f,\pi}I_{\fH_\pi}$ for some scalar $\alp_{f,\pi}$,
for each $\pi\iin\what{G}$.  Thus by the left extreme case of
(\ref{eq:hstrace}) we have
\[
\sum_{\pi\in\what{G}}d_\pi^{3/2}\norm{\hat{f}(\pi)}_2
=\sum_{\pi\in\what{G}}d_\pi\norm{\hat{f}(\pi)}_1.
\]
Hence it follows that $f\in\falg\cap\mathrm{Z}\bloneg$
if and only if $f\in\fdelg\cap\mathrm{Z}\bloneg$, with
$\fnorm{f}=\fdelnorm{f}$.  \endpf

\begin{fdelgchar2}\label{cor:fdelgchar2}
We have that $\fdelg=\falg$ if and only if $G$ admits
an abelian subgroup of finite index.
\end{fdelgchar2}

\proof  We invoke the well-known result of \cite{moore}
that
\[
\begin{matrix}G\text{ admits an abelian subgroup} \\ \text{of finite index }
\end{matrix}
\quad\iff\quad d_G=\sup_{\pi\in\what{G}}d_\pi<\infty.
\]
Now if $d_G<\infty$, then for any $u\iin\falg$ we have, using 
(\ref{eq:hstrace}), that
\[
\sum_{\pi\in\what{G}}d_\pi^{3/2}\norm{\hat{u}(\pi)}_2
\leq d_G^{1/2}
\sum_{\pi\in\what{G}}d_\pi\norm{\hat{u}(\pi)}_1
<\infty
\]
so $u\in\fdelg$.  Conversely, if $\falg=\fdelg$, then since
$\fnorm{\cdot}\leq\fdelnorm{\cdot}$, the open mapping theorem
provides us with a constant $K$ such that 
$\fdelnorm{\cdot}\leq K\fnorm{\cdot}$.
For any $\pi\iin\what{G}$ we let $\pi_{11}$ be as in (\ref{eq:piij}).
We have
\[
\hat{\pi}_{11}(\sig)=\int_G\inprod{\pi(s)\xi_1^\pi}{\xi_1^\pi}\bar{\sig}(s)ds
=\begin{cases} \frac{1}{d_\pi}\inprod{\cdot}{\bar{\xi}^\pi_1}
\bar{\xi}^\pi_1\iif\sig=\pi \\
0\text{ otherwise.}\end{cases}
\]
Thus by the rank 1 case of (\ref{eq:hstrace}) we have
\[
d_\pi^{3/2}\norm{\hat{\pi}_{11}(\pi)}_1=
d_\pi^{3/2}\norm{\hat{\pi}_{11}(\pi)}_2=\fdelnorm{\hat{\pi}_{11}}
\leq K\fnorm{\hat{\pi}_{11}}=Kd_\pi\norm{\hat{\pi}_{11}(\pi)}_1
\]
so $d_\pi\leq K^2$.  Thus $d_G\leq K^2<\infty$. \endpf

We remark that despite the identification $(G\cross G)/\Del\cong G$,
the result above tells us that for the Fourier algebra
over this coset space, $\fal{G\cross G/\Del}\cong\fdelg$
is not naturally ismorphic to $\falg$.
We will see in Section \ref{sec:weakamen} that 
$\fdelg$ can fail to be operator weakly amenable, while
$\falg$ is always operator weakly amenable.

\subsection{Some subalgebras of $\fdelg$} Let us begin
with a variant of Proposition \ref{prop:optensprod2}.

\begin{optensprod}\label{prop:optensprod}
{\bf (i)} There is a completely isometric identification 
\[
\fdelg\what{\otimes}\fdelh
\cong\fdelgh
\] 
given on elementary tensors by $f\otimes g\mapsto f\cross g$.

{\bf (ii)} If $K$ is a closed subgroup of $G$, and
$L$ is a closed subgroup of $H$, then we obtain
a completely isometric identification
\[
\fdel{G/K}\what{\otimes}\fdel{H/L}\cong\fdel{G\cross H/K\cross L}. 
\]
\end{optensprod}

\proof {\bf (i)} We have, using Proposition \ref{prop:optensprod2}, completely
isometric identifications
\begin{align*}
\fdelg\what{\otimes}\fdelh
&\cong\fal{G\cross G/\Del_G}\what{\otimes}\fal{H\cross H/\Del_H} \\
&\cong\fal{G\cross G\cross H\cross H/\Del_G\cross \Del_H}.
\end{align*}
And we have a completely isometric identification
\[
\fdelgh\cong\fal{G\cross H\cross G\cross H/\Del_{G\times H}}.
\]
Thus we must show that
\begin{equation}\label{eq:identity}
\fal{G\cross G\cross H\cross H/\Del_G\cross \Del_H}
\cong \fal{G\cross H\cross G\cross H/\Del_{G\times H}}.
\end{equation}
Let $\varsigma$ be the topological group isomorphism
$(s_1,t_1,s_2,t_2)\mapsto(s_1,s_2,t_1,t_2):
G\cross H\cross G\cross H\to G\cross G\cross H\cross H$.
This map induces a completely isometric 
isomorphism $u\mapsto u\comp\varsigma$ from
$\fal{G\cross G\cross H\cross H}$ to $\fal{G\cross H\cross G\cross H}$.
Moreover, the following digram commutes.
\[
\xymatrix{
\fal{G\cross G\cross H\cross H} \ar@<.5ex>[rr]^{u\mapsto u\comp\varsigma}
\ar[d]_{P_{\Del_G\times\Del_H}}
& & \fal{G\cross H\cross G\cross H} \ar[d]^{P_{\Del_{G\times H}}} 
\ar@<.5ex>[ll]^{u\mapsto u\comp\varsigma^{-1}}  \\
\fal{G\cross G\cross H\cross H\!:\!\Del_G\cross \Del_H}
\ar@<.5ex>[rr]^{u\mapsto u\comp\varsigma}
& & \fal{G\cross H\cross G\cross H\!:\!\Del_{G\times H}}  
\ar@<.5ex>[ll]^{u\mapsto u\comp\varsigma^{-1}}  }
\]
Since each of the top row maps are complete isometries, and
each of the maps $P_{\Del_G\times\Del_H}$ and $P_{\Del_{G\times H}}$
are complete quotient maps, we have that the bottom row maps must
each be complete quotient maps which are mutual inverses, hence
complete isometries. By standard identifications, (\ref{eq:identity}) follows.

{\bf (ii)} This can be proved exactly as Proposition \ref{prop:optensprod2}
using (i), above, in place of the identification $\falg\what{\otimes}\fal{H}
\cong\fal{G\cross H}$.  \endpf



We now define a sequence of subalgebras of $\falg$.

\begin{fdelngchar}\label{theo:fdelngchar}
Let $\fdelbg{1}=\fdelg$ and for each $n\geq 1$ let
\[
\fdelbg{n+1}=\Gam(\fdelngg).
\]

{\bf (i)} Each $\fdelng$ is a subalgebra of $\falg$, which is closed under
both left and right translations and for which all translations are
complete isometries.  Also, for each $n$, the map $\Gam:\fdelngg\to
\fdelbg{n+1}$ is a complete quotient map.

{\bf (ii)} If $f\in\fC(G)$, then
\[
f\in\fdelng\quad\iff\quad \sum_{\pi\in\what{G}}d_\pi^{(2^n+1)/2}
\norm{\hat{f}(\pi)}_2
<\infty.
\]
Moreover, the latter quantity is the norm, $\fdelnnorm{f}$.
\end{fdelngchar}

\proof  {\bf (i)} We use induction.  If the ascribed properties hold
for $\fdelng$, and we have that
\begin{equation}\label{eq:fdelgtp}
\fdelng\what{\otimes}\fdelng\cong\fdelngg.
\end{equation}
then it follows from Proposition \ref{prop:basicprop} (i)
that $\fdelbg{n+1}$ is a subalgebra of $\falg$, which is 
completely isometrically
isomorphic to a completely contractively complemented subspace
of $\fdelngg$.  It follows from Proposition \ref{prop:basicprop} (ii)
that $\fdelbg{n+1}$ is closed under left and right translations.
The formula (\ref{eq:fdelgtp}) holds for $n=1$ by Proposition
\ref{prop:optensprod} (i).  Moreover, if (\ref{eq:fdelgtp}) holds for $n$,
then we can use the proof of Proposition \ref{prop:optensprod} (ii), then (i),
to see that (\ref{eq:fdelgtp}) holds for $n+1$.

{\bf (ii)} This follows exactly as the proof of Theorem
\ref{theo:fdelgchar}, using Lemma \ref{lem:tnormcoeff} in the case $p=2$,
and induction. \endpf

We note that a simple modification of Corollary \ref{cor:fdelgchar2}
shows that if $G$ admits no abelian subgroup of finite index then
$\{\fdelng\}_{n=1}^\infty$ is a properly nested sequence
of subalgebras of $\falg$.  We also note that it follows from Proposition 
\ref{prop:spectrum} (i), that $\Sig_{\fdelg}\cong G$.  We suspect
the same holds for each $\fdelng$ ($n\geq 2$), but have not
been able to prove it.  
We cannot copy the method of proof of \cite[Thm.\ 4.1]{forrest} which
we cited in Proposition \ref{prop:spectrum} (i) to obtain that
$\Sig_{\fdelsqg}\cong G$.  Indeed, we have that
$\Del$ is not a set of spectral synthesis for $\fdelgg$ when $G$ is a nonabelian 
connected Lie group, by (\ref{eq:ggideals}), below.



\section{Amenability properties}\label{sec:weakamen}

\subsection{Failure of operator weak amenability of $\fdelg$}
We adapt arguments from
\cite[Thm.\ 7.2 \& Cor.\ 7.3]{johnson}.  We also use the
Fourier series of any $f\iin\fT(G)$ ($\fT(G)$ is defined in 
(\ref{eq:trigfunc})):
\[
f(s)=\sum_{\pi\in\what{G}}d_\pi\sum_{i,j=1}^{d_\pi}\hat{f}(\pi)_{ij}\pi_{ij}(s)
=\sum_{\pi\in\what{G}}d_\pi\mathrm{trace}\left[\hat{f}(\pi)\pi(s)^\t\right]
\]
where $\hat{f}(\pi)_{ij}=\inprod{\hat{f}(\pi)\xi^\pi_j}{\xi^\pi_i}$ and
$A^\t$ is the transpose of a matrix $A$.
This is a variant of the formula given in \cite[(34.1)]{hewittrII}.

\begin{comweakamen}\label{theo:comweakamen}
If $G$ is a nonabelian connected compact Lie group, then $\fdelg\cong
\fal{G\cross G/\Del}$ is not operator weakly amenable.
\end{comweakamen}

\proof  We will establish that $\fdelsqg$ admits a point derivation at $e$.
This implies that $\wbar{\ideal_{\fdelsqg}(\{e\})^2}
\not=\ideal_{\fdelsqg}(\{e\})$.  Indeed, the point derivation
will vanish on $\ideal_{\fdelsqg}(\{e\})^2$, but not on 
$\ideal_{\fdelsqg}(\{e\})$. Hence it follows from the construction
of $\fdelsqg$ (Theorem \ref{theo:fdelngchar}), and Corollary
\ref{cor:idealr1} (iii), that
\begin{equation}\label{eq:ggideals}
\wbar{\ideal_{\fdelggs}(\Del)^2}\not=\ideal_{\fdelggs}(\Del).
\end{equation}
This condition
implies that
$\fdelg$ is not operator weakly amenable 
by \cite[Thm.\ 3.2]{groenbaek}, which was shown to be a valid
characterisation of operator weak amenability in \cite{spronk}.

It has been shown in \cite{plymen} that under the assumptions given,
there is a subgroup $T\cong\Tee$ in $G$ such that
for any $\pi\iin\what{G}$
\[
\pi|_T\cong\bigoplus_{k=1}^{d_\pi}\chi_{n_k}\quad 
\text{with each } |n_k|< d_\pi
\]
where $\chi_n(z)=z^n$ for $n\iin\Zee\cong\what{\Tee}$.
We let $\theta\mapsto t_\theta:(-\pi,\pi]\to T$ be the parameterisation
of $T$ which corresponds to $\theta\mapsto e^{i\theta}:(-\pi,\pi]\to\Tee$.
For each $\pi\iin\what{G}$ we can choose an orthonormal basis for $\fH_\pi$
with respect to which $\pi(t_\theta)=
\mathrm{diag}(e^{in_1\theta},\dots,e^{in_{d_\pi}\theta})$;
and it follows by elementary estimates that
\[
\frac{\norm{\pi(t_\theta)^\t-I}}{\theta}
=\max_{k=1,\dots,d_\pi}
\frac{|e^{in_k\theta}-1|}{|\theta|}<d_\pi
\]
for $\theta$ in a neighbourhood of $0$.

We have for $f\iin\fT(G)$ that
\[
\frac{f(t_\theta)-f(e)}{\theta}
=\sum_{\pi\in\what{G}}d_\pi\mathrm{trace}
\left[\hat{f}(\pi)\frac{(\pi(t_\theta)^\t-I)}{\theta}\right].
\]
We note that for small $\theta$ we have, using (\ref{eq:hstrace})
\[
\left|\frac{f(t_\theta)-f(e)}{\theta}\right|
\leq\sum_{\pi\in\what{G}}d_\pi\norm{\hat{f}(\pi)}_1 
\frac{\norm{\pi(t_\theta)^\t-I}}{|\theta|}
\leq\sum_{\pi\in\what{G}}d_\pi^{5/2}\norm{\hat{f}(\pi)}_2
\]
where the last quantity is $\fdelsqnorm{f}$, 
by Theorem \ref{theo:fdelngchar} (ii).
Hence, since each $\lim_{n\to\infty}(\pi(t_\theta)-I)/\theta$ 
exists, we have that
\[
d(f)=\frac{d}{d\theta}\Big|_{\theta=0}f(t_\theta)
=\lim_{\theta\to 0}\frac{f(t_\theta)-f(e)}{\theta}.
\]
exists, and $|d(f)|\leq\fdelsqnorm{f}$ for $f\iin\fT(G)$.  
Hence $d$ extends to a contractive
point derivation on $\fdelsqg$ at $e$. \endpf

We remark that Johnson \cite[Cor. 7.3]{johnson} showed that the 
point derivation $d:\fT(G)\to\Cee$ extends to a bounded
map on $\fgamg$ (see Section \ref{sec:fdelcg}), and hence established
that $\falg$ is not weakly amenable.  Since $\Del$ is a set of spectral synthesis
for $\falgg$, we can proceed as in the first paragraph of the proof 
above to see that $d$ cannot be extended to $\fdelg$.

As an application of the above result, we obtain a new set of non-synthesis.

\begin{comweakamen2}\label{cor:comweakamen2}
Let $G$ be a compact connected non-abelian Lie group.  Then
$(\Del_G\cross\Del_G)\Del_{G\times G}$ fails spectral synthesis
for $\fal{G\cross G\cross G\cross G}$.
\end{comweakamen2}

\proof It follows from (\ref{eq:ggideals}) 
that $\Del_G$ is not a 
set of spectral synthesis for $\fdelgg$.
We then appeal to Proposition \ref{prop:optensprod} (ii)
to see that
\[
\fdelg\hat{\otimes}\fdelg\cong\fdelgg
\cong\fal{G\cross G\cross G\cross G/\Del_{G\times G}}.
\]
In the identification (\ref{eq:ggmoddel}) we have
$\Del_G\cong\{(s,s,e,e)\Del_{G\times G}:s\in G\}$.
It then follows from Corollary \ref{cor:idealr1} (iii) that 
\[
\Del_G^*\cong(\Del_G\cross\{1\})\Del_{G\times G}
=(\Del_G\cross\Del_G)\Del_{G\times G}
\]
is not spectral for $\fal{G\cross G\cross G\cross G}$. \endpf


\subsection{Failure of operator weak amenability of $\fal{G/K}$}
Let us now turn our attention to general locally compact groups.
We will let $G$ denote a locally compact group. 
Let us collect some useful facts.

\begin{restricth}\label{lem:restricth}
Let $H$ be a closed subgroup of $G$ and $K$ be a compact
subgroup of $H$.  Then the restriction map $u\mapsto u|_H$
maps $\fal{G\!:\! K}$ onto $\fal{H\!:\! K}$.
\end{restricth}

\proof This is \cite[Lem.\ 3.6 (ii)]{ghahramanis}.  \endpf


\begin{fident}\label{lem:fident}
Let $N$ be a compact normal subgroup
of $G$, and $q:G\to G/N$ be the canonical quotient map.
If $\til{K}$ is a compact subgroup of $G/N$ and $K=q^{-1}(\til{K})$,
then the algebras $\fal{G\!:\!K}$ and $\fal{G/N\!:\!\til{K}}$
are completely isometrically isomorphic.
\end{fident}

\proof The map $u\mapsto u\comp q:\fal{G/N}\to\falg$
is a complete isometry with range $\fal{G\!:\!N}$.
If $u\in\fal{G/N\!:\!\til{K}}$, then for $s\iin G$ and $k\iin K$ we have
\[
u\comp q(sk)=u\bigl(q(s)q(k)\bigr)=u\bigl(q(s)\bigr)=u\comp q(s)
\]
so $u\comp q\in\fal{G\!:\!K}$.  Conversely, let
$v\in\fal{G\!:\!K}\subset\fal{G\!:\!N}$.  Let $u\iin\fal{G/N}$ be
such that $v=u\comp q$.  For any $\til{s}\iin G/N$ and
$\til{k}\in \til{K}$, find $s\in G$ and $k\in K$ so
$q(s)=\til{s}$ and $q(k)=\til{k}$. We have
\[
u(\til{s}\til{k})=u\comp q(sk)=v(sk)=v(s)=u(\til{s})
\]
so $u\in\fal{G/N\!:\!\til{K}}$.  \endpf

We now obtain a generalisation of Theorem \ref{theo:comweakamen}.

\begin{weakamen}\label{theo:weakamen}
Suppose $G$ contains a connected
nonabelian compact subgroup $K$.  Then there is a compact subgroup
$K^*$ of $G\cross G$ such that $\fal{G\cross G/K^*}$ is not operator
weakly amenable.
\end{weakamen}

\proof There exists a closed normal subgroup $N$ of $K$
such that $K/N$ is a compact Lie group \cite[(28.61)(c)]{hewittrII}.
Moreover, we may have that $K/N$ is nonabelian.
Indeed, if $st\not=ts$ in $K$, find a neighbourhood
$U$ of $e$ in $K$ such that $ts\not\in stU$; find in $U$
a compact subgroup $N$ so $K/N$ is a Lie group.
Then by Theorem \ref{theo:comweakamen}, 
$\fal{K/N\cross K/N\!:\!\Del_{K/N}}$ is not operator weakly amenable.
Let $q:K\cross K\to K/N\cross K/N$ be the quotient map and
$K^*=q^{-1}(\Del_{K/N})$.  Then by Lemma \ref{lem:fident},
$\fal{K\cross K\!:\!K^*}\cong\fal{K/N\cross K/N\!:\!\Del_{K/N}}$
is not operator weakly amenable.  Moreover, by Lemma
\ref{lem:restricth}, $\fal{G\cross G/K^*}\cong \fal{G\cross G\!:\!K^*}$ admits
$\fal{K\cross K\!:\!K^*}$ as a quotient, and hence is not
operator weakly amenable either.  \endpf


\subsection{Examples of hyper-Tauberian $\fal{G/K}$}\label{ssec:hypertaub}
In this section 
we will always let $G$ denote a locally compact group.
In this section we generalise the fact that $\falg$ is hyper-Tauberian
when the connected component of the identity
is abelian \cite[Theo.\ 21]{samei2}.  
Our approach is inspired by that of
\cite[Thm.\ 3.3]{forrestr}.  However, in dealing with coset spaces
some extra technicalities arise.
The following lemma deals with
some of these technicalities.

\begin{weird}\label{lem:weird}
Suppose $G$ contains an open subgroup $G_0$ and
a compact subgroup $K$ for which 
$\fal{G_0/G_0\cap K}$ is hyper-Tauberian.
Then $\fal{G/K}$ is hyper-Tauberian.
\end{weird}

\proof We will identify $\fal{G/K}\cong\fal{G\!:\! K}$, etc., so we
may work within the algebra $\falg$.

Let $H=G_0\cap K$.  We will first establish that
$\fal{G_0\!:\! H}$ is boundedly isomorphic
to a certain subalgebra of $\fal{G\!:\! K}$.  

Since $H$ is open in $K$, it is of finite index.
Thus there is a finite set  $\{k_1,k_2,\dots,k_n\}\subset K$
for which
\begin{equation}\label{eq:hcoset}
K=\bigcup_{i=1}^nHk_i\aand Hk_i\cap Hk_j=\varnothing\iif i\not=j.
\end{equation}
It then follows easily that
\[
G_0K=\bigcup_{i=1}^nG_0k_i\aand G_0k_i\cap G_0k_j=\varnothing\iif i\not=j
\]
and thus $GK$ is a union of open cosets.  Then
the indicator function $1_{GK}\in\fsalg$, and it is clear that
$k\mult 1_{G_0K}=1_{G_0K}$ for each $k\in K$.  Hence
$1_{G_0K}\fal{G\!:\!K}$ is a closed subalgebra of $\fal{G\!:\!K}$.

Now if $u\in\fal{G\!:\! K}$, then $u|_{G_0}\in\fal{G_0\!:\! H}$
by Lemma \ref{lem:restricth}.  The restriction map 
\begin{equation}\label{eq:restrict}
u\mapsto u|_{G_0}:1_{G_0K}\fal{G\!:\! K}\to\fal{G_0\!:\! H}
\end{equation}
is injective, for if $u\iin\fal{G\!:\! K}$ satisfies
$u(s)\not=0$, for some $s=s_0k\iin G_0K$, where
$s_0\in G_0$ and $k\in K$, then $u(s_0)=u(s)\not=0$.
Let us see the map in (\ref{eq:restrict}) is surjective.  Indeed,
if $v\in\fal{G_0\!:\! H}$ we define $\til{v},w\iin\falg$ by
\[
\til{v}(s)=\begin{cases} v(s) &\iif s\in G_0 \\ 0 & \text{ otherwise}
\end{cases} \quad\aand\quad w=\sum_{i=1}^nk_i^{-1}\mult\til{v}.
\]
Then $w\in1_{G_0K}\falg$ with $w|_{G_0}=v$.
Now, if $s_0\in G_0$, $h\in H$ and $k_i$ is as above, then,
since $s_0h\in G_0$ and $v\in\fal{G_0\!:\!H}$, we have
\[
w(s_0hk_i)=v(s_0h)=v(s_0)=w(s_0)
\]
and hence by (\ref{eq:hcoset}), $w(s_0k)=w(s_0)$ for $s_0\iin G_0$ 
and $k\iin K$.  Thus if $s\in G$, and $k\in K$, then either $s\in G_0K$
and we can find
$s_0\iin G_0$ and $k_i$ as above such that $s=s_0k_i$, so we have
\[
w(sk)=w(s_0k_ik)=w(s_0)=w(s_0k_i)=w(s);
\]
or $s\not\in G_0K$, so for any $k\iin K$, $sk\not\in G_0K$
and thus $w(sk)=0=w(s)$.  We conclude that
$w\in\fal{G\!:\! K}$.  Thus the map in (\ref{eq:restrict}) 
is a contractive bijection
which is also a homomorphism.  It follows from the open mapping
theorem that this map is an isomorphism. Thus,
since $\fal{G_0\!:\!H}$ is hyper-Tauberian, $1_{G_0K}\fal{G\!:\! K}$ 
is hyper-Tauberian.

Now let $T:\fal{G\!:\!K}\to\fal{G\!:\!K}^*$ be a bounded local
operator, i.e.\ so $N^*\comp T\comp N:\fal{G/K}\to\fal{G/K}^*$
is a local operator.  We will show for
$u_1,u_2,u_3\iin\fal{G\!:\!K}$, each having compact support, that
\begin{equation}\label{eq:tuvw}
\dpair{T(u_1u_2)}{u_3}=\dpair{u_1T(u_2)}{u_3}.
\end{equation}
Then it follows from Proposition \ref{prop:spectrum} (ii), that
(\ref{eq:tuvw}) holds for any $u_1,u_2,u_3$ in $\fal{G\!:\!K}$.

Since $G_0$ is open, there are $t_1,\dots,t_n\iin G$ such that 
\[
\bigcup_{j=1}^3\supp{u_j}
\subset\bigcup_{i=1}^n t_iG_0\subset \bigcup_{i=1}^n t_iG_0K.
\]
The map $u\mapsto t_i\con u$ is an isometric isomorphism
from $1_{G_0K}\fal{G\!:\! K}$ to $1_{t_iG_0K}\fal{G\!:\! K}$
for each $i$, so each $1_{t_iG_0K}\fal{G\!:\! K}$ is hyper-Tauberian.
Now let 
\[
w_1=1_{t_1G_0K},\aand w_i=1_{t_iG_0K}\left(1-\sum_{k=1}^{i-1}w_k\right)
\ffor i=2,\dots,n.
\]
Then each $w_i$ is an idempotent in $\fsalg$ with $k\mult w_i=w_i$
for each $k\iin K$.  Moreover, $w_iw_j=0$ if $i\not=j$.  For each $i$
the map $u\mapsto w_iu$ from $1_{t_iG_0K}\fal{G\!:\! K}
\to w_i\fal{G\!:\! K}$ is a surjective homomorphism, so
$w_i\fal{G\!:\! K}$ is hyper-Tauberian by \cite[Thm.\ 12]{samei2}.  
Then the algebra
\[
\fA=\sum_{i=1}^nw_i\fal{G\!:\! K}
\cong\bigoplus_{i=1}^nw_i\fal{G\!:\! K}
\]
is hyper-Tauberian, by \cite[Cor.\ 13]{samei2}, and contains each $u_j$.
The inclusion map $\iota:\fA\to\fal{G\!:\!K}$ is an $\fA$-module
map, so $\iota^*\comp T\comp\iota:\fA\to\fA^*$ is an $\fA$-local
operator.  Hence $\iota^*\comp T\comp\iota$ is
an $\fA$-module map and (\ref{eq:tuvw}) holds.  \endpf

\begin{weakamens}\label{theo:weakamens}
If $G$ has abelian connected component of the identity, $G_e$,
then for any compact subgroup $K$, $\fal{G/K}$ is hyper-Tauberian.
\end{weakamens}

\proof  As in the proof of the lemma above, we will identify
$\fal{G/K}\cong\fal{G\!:\!K}$, etc.

We will first assume that $G$ is almost connected.

Let $U$ be a neighbourhood identity in $G$. Then, 
by \cite{montgomeryz}, there
is a compact normal subgroup $N_U\subset U$ such that
$\til{G}=G/N_U$ is a Lie group.  Let $q:G\to\til{G}$ be the canonical
quotient map and $\til{K}=q(K)$, so $\til{K}$ is a compact 
subgroup of $\til{G}$.  We have that the connected component
of the identity of $\til{G}$ satisfies $\til{G}_e=\wbar{q(G_e)}$
by \cite[(7.12)]{hewittrI}.
Then $\til{G}_e$ is abelian, and open in $\til{G}$ since the latter
is a Lie group.  It then follows that $\til{G}_e\cap\til{K}$ is
normal in $\til{G}_e$ and $\fal{\til{G}_e\!:\!\til{G}_e\cap\til{K}}
\cong\fal{\til{G}_e/\til{G}_e\cap\til{K}}$ is hyper-Tauberian
by \cite[Prop.\ 18]{samei2}.
Then it follows from Lemma \ref{lem:weird} that
$\fal{\til{G}\!:\!\til{K}}$ is hyper-Tauberian.

We let $K_U=q^{-1}(\til{K})=KN_U$.
Then Lemma \ref{lem:fident} tells us that
$\fal{G\!:\! K_U}$ is isometrically isomorphic
to $\fal{\til{G}\!:\!\til{K}}$, and hence is hyper-Tauberian.
Since $K_U\supset K$, it follows that $\fal{G\!:\!K_U}\subset
\fal{G\!:\! K}$.

Let $u\in\fal{G\!:\! K}$, and $\eps>0$.  
Fix a compact neighbourhood $V$ of $e$.  Find a 
neighbourhood $U$ of $e$ such that $U\subset V$ and
\[
\fnorm{s\mult u-u}<\frac{\eps}{m(KV)}\ffor s\in U.
\]
Then find normal a subgroup $N_U\subset U$ as above and let
\[
P_Uu=\int_{K_U}k\mult u\,dk\in\fal{G\!:\!K_U}.
\]
We note that for any $k=k'n\in KN_U$ with $k'\iin K$ and
$n\iin N_U$ we have
\[
\fnorm{k\mult u-u}=\fnorm{k'\mult (n\mult u)-u}
=\fnorm{n\mult u-k'^{-1}\mult u}=\fnorm{n\mult u-u}<\frac{\eps}{m(KV)}
\]
since $u\in\fal{G\!:\! K}$ and right translation on $\falg$ is an isometry.  
Hence we find
\[
\fnorm{P_Uu-u}\leq\int_{K_U}\fnorm{k\mult u-u}dk< \frac{\eps}{m(KV)}m(K_U)
\leq \eps.
\]
Thus we have $\lim_{U\searrow \{e\}}P_Uu=u$.

Now let $T:\fal{G\!:\!K}\to\fal{G\!:\!K}^*$ be a bounded local operator.
Since the inclusion $\iota: \fal{G\!:\!K_U}\to\fal{G\!:\!K}$
is an $\fal{G\!:\!K_U}$-module map, and $\fal{G\!:\!K_U}$
is hyper-Tauberian, $\iota^*\comp T\comp\iota:
\fal{G\!:\!K_U}\to\fal{G\!:\!K_U}^*$ is an $\fal{G\!:\!K_U}$-local
map and hence an $\fal{G\!:\!K_U}$-module map.
Hence if $u_1,u_2,u_3\in \fal{G\!:\!K}$, then
\[
\dpair{T(P_Uu_1\,P_Uu_2)-P_Uu_1T(P_Uu_2)}{P_Uu_3}
=0
\]
Taking $U\searrow\{e\}$, as above, we obtain
\[
\dpair{T(u_1u_2)-u_1T(u_2)}{u_3}=0
\]
and hence $T$ is an $\fal{G\!:\!K}$-module map too. 

Finally, if $G$ is not almost connected, we can find an almost
connected open subgroup $G_0$.
Then, from above, $\fal{G_0\!:\! G_0\cap K}$ is hyper-Tauberian.
We then appeal to Lemma \ref{lem:weird}.
\endpf

We say $G$ is a [SIN]-group if there is a neighbourhood basis
at $e$ consisting of sets invariant under inner automorphisms.
We obtain for such groups, a partial converse of Theorem \ref{theo:weakamens},
which is similar to \cite[Theo.\ 3.7]{forrestr}.

\begin{weakamens2}
If $G$ is a [SIN]-group, then $\fal{G\cross G/K^*}$ is hyper-Tauberian
for every compact subgroup $K^*$ of $G\cross G$ if and only if
$G_e$ is abelian.
\end{weakamens2}

\proof Sufficiency is an obvious consequence of Theorem 
\ref{theo:weakamens}.  To see necessity, we first note
that $G_e$ is a [SIN]-group and the Freudenthal-Weil Theorem
(see \cite[12.4.28]{palmerII}) tells us that $G_e\cong V\cross K$
where $V$ is a vector group and $K$ a connected compact group.
If $K$ is nonabelian, we appeal to Theorem \ref{theo:weakamen}
to see that there exists a subgroup $K^*$ of $G\cross G$
such that $\fal{G\cross G/K^*}$ is not operator weakly amenable,
hence not weakly amenable and not hyper-Tauberian.  \endpf

\subsection{Operator amenability of $\fal{G/K}$}\label{ssec:opamen}
We recall that [MAP]$_K$-groups were defined in Section \ref{ssec:idealr}.
We note that $\fT_f(G)\falg\subset\falg$, since $\fT_f(G)$ is
a subalgebra of the Fourier-Stieltjes algebra $\mathrm{B}(G)$.

\begin{opamen}\label{theo:opamen}
Let $G$ be an amenable locally compact group and $K$ a compact subgroup
so that $G$ is [MAP]$_K$.  Then
the following are equivalent:

{\bf (i)} $\fal{G/K}$ is operator amenable; and

{\bf (ii)} $(K\cross K)\Del$ is in the closed coset ring of $G\cross G$.

\noindent Moreover, If $G$ is compact, each of the above is equivalent to

{\bf (iii)} $\ideal_{\fdelg}(K)$ has a bounded approximate identity.
\end{opamen}

\proof Since the map $u\mapsto\check{u}$ is an
isomorphism on $\falgg$, $\ideal_{\fal{G\times G}}((K\cross K)\Del)$
has a b.a.i.\ if and only if $\ideal_{\fal{G\times G}}(\Del(K\cross K))$
has a b.a.i.  It follows from \cite{forrestkls} that 
$\ideal_{\fal{G\times G}}(\Del(K\cross K))$ has a b.a.i.\ if and only
if $\Del(K\cross K)$ is in the closed coset ring of $G\cross G$.
An application of Corollary \ref{cor:idealr1} (i),
tells us that $\ideal_{\fal{G\times G}}(\Del(K\cross K))$ has a b.a.i.\
if and only if 
$\ideal_{\fal{G\times G/K\times K}}(\Del_{G/K})$
has a b.a.i., where $\Del_{G/K}=\{(s,s)K\cross K:s\in G\}$.
Under the identification
\[
\fal{G/K}\what{\otimes}\fal{G/K}\cong\fal{G\cross G/K\cross K}
\]
provided by Proposition \ref{prop:optensprod2}, 
$\ideal_{\fal{G\times G/K\times K}}(\Del_{G/K})$
corresponds to the kernel of the multiplication map
$m:\fal{G/K}\what{\otimes}\fal{G/K}\to\fal{G/K}$.  Since
$\fal{G/K}$ has a bounded approximate identity, this kernel
has a bounded approximate identity if and only if $\fal{G/K}$ is
operator amenable, by the completely bounded version of
a splitting result from \cite{helemski} (also see \cite[3.10]{curtisl}).  
Hence the equivalence of (i) and (ii), above, is established.

If $G$ is compact then by Corollary \ref{cor:idealr1} (i), 
$\ideal_{\fdelg}(K)$ has a b.a.i.\
if and only if $\ideal_{\fal{G\times G}}(K^*)$ has a b.a.i.,
where
\begin{equation}\label{eq:kstar}
K^*=(K\cross\{e\})\Del=(K\cross K)\Del.
\end{equation}
Hence it again follows \cite{forrestkls} that properties (ii) and (iii), 
above, are equivalent.   \endpf

The situation that $(K\cross K)\Del$ is in the coset ring of $G\cross G$
seems rare.  It does occur, for example, when $K$ contains
a subgroup $N$, of finite index, which is normal in $G$.  Thus it
is only in such cases that we know $\fal{G/K}$ is operator amenable.

However, we gain some situations
in which $(K\cross K)\Del$ is a set of local synthesis
for $\falgg$.

\begin{specsynth}
If $G$ has abelian connected component of the identity and 
$K$ is a compact subgroup of $G$ so that $G$ is [MAP]$_K$, then 

{\bf (i)} $(K\cross K)\Del$ is a set of local synthesis for $\falgg$.

\noindent Moreover, if $G$ is compact then

{\bf (ii)} $K$ is a set of spectral synthesis for $\fdelg$.
\end{specsynth}

\proof {\bf (i)}
By Theorem \ref{theo:weakamens}, $\fal{G/K}$ is hyper-Tauberian.  Hence
it is operator hyper-Tauberian.  Thus, by \cite[Thm.\ 6]{samei2},
$\Del_{G/K}=\{(s,s)K\cross K:s\in G\}$, is a set of local synthesis
for $\fal{G\cross G/K\cross K}$, since the latter 
is isomorphic to $\fal{G/K}\what{\otimes}\fal{G/K}$ by Proposition
\ref{prop:optensprod2}.  (We note that \cite[Thm.\ 6]{samei2} is proved
for the projective tensor product of a hyper-Tauberian algebra with itself.
However, an inspection of the proof, coupled with the formula representing
an arbitrary element of the operator projective tensor product
in \cite[10.2.1]{effrosrB},
shows that it holds for an operator hyper-Tauberian algebra with the operator
projective tensor product.)
Then it follows from Corollary \ref{cor:idealr1}
(ii') that $\Del_{G/K}^*=\Del(K\cross K)$ is a set of local synthesis
for $\falgg$.  Since $u\mapsto\check{u}$ is an isomorphism on
$\falgg$, $(K\cross K)\Del$ is also
a set of local synthesis.  

{\bf (ii)}  By Corollary \ref{cor:idealr1} (ii), $K$ is spectral
for $\fdelg\cong\fal{G\cross G/\Del}$ if and only if
$K^*=(K\cross K)\Del$ (see (\ref{eq:kstar})) is spectral for
$\falgg$.  We appeal to (i), above.  \endpf


\section{Convolution}\label{sec:fdelcg}

\subsection{Convolution on compact groups}
We close this article by addressing, in part, the case of what happens if we
replace the map $\Gam$, in Section \ref{ssec:twisted}, with convolution.   

We let $\fA(G\cross G)$ be as in Section \ref{ssec:twisted}, 
and insist further that
the group action of left translation  is isometric
on $\fA(G\cross G)$ and continuous on $G\cross G$.  We then define
a group action $(r,f)\mapsto r\diamond f:
G\cross \fA(G\cross G)\to\fA(G\cross G)$ by
\[
r\diamond w(s,t)=w(sr,r^{-1}t)=(r,e)\mult[(e,r)\con w](s,t)
=(e,r)\con[(r,e)\mult w](s,t).
\]
We let $\check{\Del}=\{(t,t^{-1}):t\in G\}$ and define
\[
\fA(G\cross G\!:\!\check{\Del})=\{f\in\fA(G\cross G):r\diamond f=f
\text{ for every }r\iin G\}.
\]
We note that for $w\in\fA(G\cross G\!:\!\check{\Del})$,
$w(s,t)=w(s_1,t_1)$, provided $(ss_1^{-1},tt_1^{-1})\in\check{\Del}$,
even though $\check{\Del}$ is not a subgroup unless $G$ is abelian.
We then define
\begin{align*}
\check{P}:\fA(G\cross G)\to\fA(G\cross G)&,\qquad
\check{P}w=\int_Gr\diamond wdr \\ \aand
\check{M}:\fA(G\cross G\!:\!\check{\Del})\to\fC(G)&,\qquad
\check{M}f(s)=f(s,e).
\end{align*}
Then $\check{P}$ is a contractive idempotent whose range is
$\fA(G\cross G\!:\!\check{\Del})$, in particular $\check{P}$ is a quotient
map.  The map $\check{M}$ is injective; we denote its range
$\fA_{\check{\Del}}(G)$ and assign it the norm which
makes $\check{M}$ an isometry.  Then $\check{M}$ has inverse
\[
\check{N}:\fA_{\check{\Del}}(G)\to\fA(G\cross G\!:\!\check{\Del}),
\qquad \check{N}f(s,t)=f(st).
\]
Finally, we define 
\[
\check{\Gam}:\fA(G\cross G)\to\fA_{\check{\Del}}(G),
\qquad \check{\Gam}=\check{M}\comp\check{P}.
\]
If $\fA(G\cross G)$ contains an elementary product $f\cross g$, then
$\check{\Gam}f\cross g=f\con g$.  Thus, $\check{\Gam}$ may be regarded simply
as the convolution map.

\subsection{Convolution on the Fourier algebras}
We will consider, now, only the case where $\fA(G\cross G)=\falgg$.
As in Section \ref{ssec:twisted}, it is easy to verify that $\check{\Gam}:
\falgg\to\fcdelg$ is a complete quotient map and
$\fcdelg\subset\falg$.

We recall that $\fgamg$ is the subalgebra of $\falg$ defined in 
Section \ref{ssec:twisted}, Example (iii).
For this algebra we have
\begin{equation}\label{eq:johnson}
f\in\fgamg\quad\iff\quad
\sum_{\pi\in\what{G}}d_\pi^2\norm{\hat{f}(\pi)}_1<\infty
\end{equation}
and the latter quantity is the norm $\fgamnorm{f}$, by 
\cite[Prop.\ 2.5]{johnson}.

\begin{convolution}
We have $\fcdelg=\fgamg$, isometrically. 
Morevoer, $\fcdelg=\falg$ if and only if
$G$ admits an abelian subgroup of finite index.
\end{convolution}

\proof 
We have computation similar to that in Lemma
\ref{lem:tnormcoeff}.  
For $f\iin\fC(G)$ and $\pi\iin\what{G}$ we have
\begin{align*}
\what{\check{N}f}(\bar{\pi}\cross\bar{\pi})
&=\int_G\int_G f(st)\pi(s)\otimes\pi(t)dt \\
&=\int_G\int_G f(s)\pi(st^{-1})\otimes\pi(t)dt \\
&=\left[\int_G f(s)\pi(s)\otimes\pi(e)ds\right]\comp
\left[\int_G\pi(t^{-1})\otimes\pi(t)dt\right] \\
&=[\hat{f}(\pi)\otimes I_{\fH_\pi}]\comp \frac{1}{d_\pi}U
\end{align*}
where $U$ is a unitary, in fact a permutation matrix,
as we shall see below.  Indeed, identify
the linear operators on $\fH_\pi$ with the matrix space
$M_{d_\pi}$ via an orthonormal basis, and then identify
$M_{d_\pi}\otimes M_{d_\pi}\cong M_{d_\pi^2}$. 
We obtain, using (\ref{eq:orthrel}), $ij,kl$th entry
\begin{align*}
\left(\int_G\pi(t^{-1})\otimes\pi(t)dt\right)_{ij,kl}
&=\int_G\pi_{ij}(t^{-1})\pi_{kl}(t)dt \\
&=\int_G\wbar{\pi_{ji}(t)}\pi_{kl}(t)dt 
=\frac{1}{d_\pi}\del_{il}\del_{jk}
\end{align*}
where $\del_{il}$ and $\del_{jk}$ are the Kronecker $\del$-symbols.

Thus it follows that
\[
\norm{\what{\check{N}f}(\bar{\pi}\cross\bar{\pi})}_1
=\frac{1}{d_\pi}\norm{\hat{f}(\pi)\otimes I_{\fH_\pi}}_1
=\norm{\hat{f}(\pi)}_1.
\]
If $\pi\not=\sig$ in $\what{G}$ then it can be shown, just as above,
that $\what{\check{N}f}(\pi\cross\sig)=0$.

We  then obtain that
\[
f\in\ran\check{\Gam}\quad\iff\quad
\sum_{(\pi,\sig)\in\what{G\times G}}d_\pi d_\sig
\norm{\what{\check{N}f}(\pi\cross\sig)}_1
=\sum_{\pi\in\what{G}}d_\pi^2\norm{\hat{f}}_1<\infty.
\]
This is precisely the characterisation obtained
for $\fgamg$ in (\ref{eq:johnson}).  

It can be shown, similarly as in Corollary \ref{cor:fdelgchar2},
that $\fgamg=\falg$ if and only if $G$ has an abelian subgroup
of finite index.  \endpf

Note that it follows from (\ref{eq:hstrace}) and Theorems \ref{theo:fdelgchar}
and \ref{theo:fdelngchar} that there are contractive inclusions
\[
\fdelsqg\subset\fgamg\subset\fdelg.
\]
Also, note that since $u\mapsto\check{u}$ is an isometric isomorphism
on $\falg$, the definition of $\fgamg$ given in Section \ref{ssec:twisted},
Example (iii), yields the equality
\[
\check{\Gam}\bigl(\falg\what{\otimes}\falg\bigr)
=\check{\Gam}\bigl(\falg\otimes^\gam\falg\bigr).
\]

{
\bibliography{convalgbib}
\bibliographystyle{plain}
}

\medskip
Address: {\sc Department of Pure Mathematics, University of Waterloo,
Waterloo, ON\quad N2L 3G1, Canada} 

\medskip
E-mail addresses:  {\tt beforres@uwaterloo.ca, esamei@uwaterloo.ca,  
\linebreak nspronk@uwaterloo.ca}

\end{document}